\documentclass[11pt]{amsproc}
\usepackage{amssymb}
\usepackage{amsfonts}

\theoremstyle{plain}

\newtheorem{corollary}{Corollary}

\newtheorem{definition}{Definition}

\newtheorem{lemma}{Lemma}

\newtheorem{proposition}{Proposition}

\newtheorem{theorem}{Theorem}
\numberwithin{equation}{section}
\input{tcilatex}

\begin{document}
\title[On the second Paneitz-Branson invariant]{On the second
Paneitz-Branson invariant}
\author{Mohammed Benalili}
\address{Departement of Mathematics, Faculty of Sciences BP119 University
Abou-Bekr Belka\"{\i}d, Tlemcen Algeria.}
\email{m\_benalili@mail.univ-tlemcen.dz}
\author{Hichem Boughazi }

\begin{abstract}
We define the second Paneitz-Branson operator on a compact Einsteinian
manifold of dimension $n\geq 5$ and we give sufficient conditions that make
it attained.
\end{abstract}

\maketitle

\section{\ Introduction}

In 1983, Paneitz \cite{11} discovered a conformally invariant fourth order
operator on $4$-dimensional Riemannian manifolds. Branson\cite{2} extended
the notion to Riemannian manifolds of dimension $n\geq 5$. This operator has
geometrical roots, it is associated to the notion of the $Q$-curvature which
can be seen as the analogue of the scalar curvature for the conformal
Laplacian. Let $(M,g)$ be a Riemannian manifold; the Paneitz-Branson
operator reads as%
\begin{equation*}
P_{g}(u)=\Delta ^{2}u-div\left( \frac{(n-2)^{2}+4}{2(n-1)(n-2)}S_{g}-\frac{4%
}{n-2}Ric_{g}\right) du+\frac{n-4}{2}Q_{g}u
\end{equation*}%
where $Ric_{g}$ and $S_{g}$ denote respectively the Ricci curvature and the
scalar curvature of $g$ and where%
\begin{equation*}
Q_{g}=\frac{1}{2(n-1)}\Delta S_{g}+\frac{n^{3}-4n^{2}+16(n-1)}{%
8(n-1)^{2}(n-2)^{2}}S_{g}^{2}-\frac{2}{(n-2)^{2}}\left\vert
Ric_{g}\right\vert ^{2}\text{.}
\end{equation*}%
The conformal property of the Paneitz-Branson expresses as: let $\widetilde{g%
}=\varphi ^{\frac{4}{n-4}}g$ be a conformal metric to $g$, where $\varphi >0$
is smooth function on $M$. Then 
\begin{equation*}
P_{g}(u\varphi )=\varphi ^{N-1}P_{g}(u)
\end{equation*}%
where $N=\frac{2n}{n-4}$.

Observe that when $(M,g)$ is Einstein, the Paneitz-Branson operator is
reduced to%
\begin{equation*}
P_{g}(u)=\Delta ^{2}u+\alpha \Delta u+\overline{\alpha }u
\end{equation*}%
where 
\begin{equation*}
\Delta =-div\nabla
\end{equation*}%
and 
\begin{equation*}
\alpha =\frac{n^{2}-2n-4}{2n(n-1)}S_{g}\text{ ,\ \ }\overline{\alpha }=\frac{%
(n-4)(n^{2}-4)}{16n(n-1)^{2}}S_{g}^{2}\text{.}
\end{equation*}%
Notice that 
\begin{equation}
\frac{\alpha ^{2}}{4}-\overline{\alpha }=\frac{S_{g}^{2}}{n^{2}(n-1)^{2}}%
\text{.}  \label{1}
\end{equation}%
Let $H_{2}^{2}(M)$ be the standard Sobolev space, which is the completion of
the space 
\begin{equation*}
C_{2}^{2}(M)=\left\{ \varphi \in C^{\infty }(M)\text{, }\left\Vert \varphi
\right\Vert _{2,2}<\infty \right\}
\end{equation*}%
with respect to the norm%
\begin{equation*}
\left\Vert \varphi \right\Vert _{2,2}=\left( \left\Vert \Delta \varphi
\right\Vert _{2}^{2}+\left\Vert \nabla \varphi \right\Vert
_{2}^{2}+\left\Vert u\right\Vert _{2}^{2}\right) ^{\frac{1}{3}}.
\end{equation*}%
Let $Gr_{k}(H_{2}^{2})$ be the $k$-dimensional Grassmannian manifold in $%
H_{2}^{2}(M)$ i.e. the set of all subspaces of $H_{2}^{2}(M)$ of dimension $%
k\geq 1$.

Denote by $\left[ g\right] $ the conformal class of the metric $g$ i.e. $%
\widetilde{g}\in \left[ g\right] $, $\widetilde{g}=ug$ with $u>0$ a smooth
function on $M.$ The minimax characterization of the eigenvalue of order $%
k\geq 1$ of the Paneitz-Branson operator $P_{g}$ is given by%
\begin{equation*}
\lambda _{k}(g)=\inf_{V\in G^{r}(H_{2}^{2})}\sup_{v\in V-\left\{ 0\right\} }%
\frac{\int_{M}vP_{g}(v)dv_{g}}{\int_{M}v^{2}dv_{g}}\text{.}
\end{equation*}%
Similarly to the Yamabe invariant of higher order introduced by Amman and
Humbert($\left[ 1\right] )$, we define the Paneitz-Branson invariant.

\begin{definition}
Let $k\in N^{\ast }$. The $k^{th}$ Paneitz-Branson invariant is defined by%
\begin{equation*}
\mu _{k}(M,g)=\inf_{\widetilde{g}\in \left[ g\right] }\lambda _{k}(%
\widetilde{g})Vol(M,\widetilde{g})^{\frac{4}{n}}.
\end{equation*}
\end{definition}

In a recent paper \cite{1} Amman and Humbert introduced the Yamabe invariant
of high order $\mu _{k}(M,g)$, $k\geq 1$ and studied $\mu _{2}(M,g)$, mainly
they showed that contrary to the standard Yamabe invariant $\mu _{1}(M,g)$
the second invariant $\mu _{2}(M,g)$ cannot be attained by a metric if the
manifold $(M,g)$ is connected. To find a minimizer to $\mu _{2}(M,g)$, they
enlarge the class $\left[ g\right] $ of conformal metric to what they called
the generalized conformal metric to $g$ \ i.e. $\widetilde{g}\in \left[ g%
\right] $ if $\widetilde{g}=u^{2^{\ast }-2}g$ where $u\in L^{2^{\ast }}(M)$
and $u\geq 0$ not indentically null and where $2^{\ast }=\frac{2n}{n-2}$.

The goal of this paper is to study the second Paneitz-Branson invariant on
Einsteinian manifolds we seek for situations where this latter is attained.
Observe that to have positive solutions in case of the Yamabe invariant it
suffices to remark that for any $u\in H_{1}^{2}(M)$, $\left\vert
u\right\vert \in H_{1}^{2}(M)$ and $\left\vert \nabla \left\vert
u\right\vert \right\vert =\left\vert \nabla u\right\vert $ which is no
longer true in the case of the Branson-Paneitz operator because of the term $%
\int_{M}\left( \Delta u\right) ^{2}dv_{g}$ and also if $u\in H_{2}^{2}(M)$, $%
\left\vert u\right\vert $ is not necessary in $H_{2}^{2}(M)$. The condition(%
\ref{1}) implies by (\cite{10} Theorem1.1) that Paneitz-Branson operator is
coercive i.e.%
\begin{equation*}
\int_{M}uP_{g}(u)dv_{g}\geq \Lambda \left\Vert u\right\Vert _{2,2}
\end{equation*}%
where the left hand side of this inequality has to be understood in the
distribution sense and where $\Lambda >0$ is a constant.

Hereafter, the space $H_{2}^{2}(M)$ will be endowed with the norm%
\begin{equation*}
\left\Vert u\right\Vert =\left( \int_{M}uP_{g}(u)dv_{g}\right) ^{\frac{1}{2}}
\end{equation*}%
which is equivalent to the norm $\left\Vert .\right\Vert _{2,2}$.

$\left\Vert ,\right\Vert _{p}$ will denote the $L^{p}$-norm with respect to
the Riemannian measure $dv_{g}$.

The main results we obtain are

\begin{theorem}
If the compact manifold $(M,g)$ is Einstein and of dimension $n\geq 12$ then 
$\mu _{2}(M,g)$ is attained by a generalized metric.
\end{theorem}

\begin{theorem}
Let $(M,g)$ be a compact Einstein manifold of positive scalar curvature and
of dimension $n\geq 5$. Assume that $\mu _{2}(M,g)$ is attained by a
generalized metric $u^{N-2}g$ with $u\in L^{N}(M)$ and $u\geq 0$ not
identically null.

Then there exist a nodal solution $w\in C^{4,\alpha }(M)$ ($\alpha <N-2$) to
the equation $P_{g}(w)=\mu _{2}(M,g)u^{N-2}w$ such that $\left\vert
w\right\vert =u$.
\end{theorem}

Our paper is organized as follows

In the first section we give some properties of the first and second
eigenvalues of the Branson-Paneitz operator. In the second one we establish
a Sobolev inequality related to the second Branson-Paneitz invariant $\mu
_{2}(M,g)$. The third section is devoted to the existence of a minimizer to $%
\mu _{2}(M,g)$. In the fourth section an estimation of $\mu _{2}(M,g)$ is
given in terms of $\mu _{1}(M,g)$ and of the best constant $K_{2}$ in the
Sobolev embedding of $H_{2}^{2}(R^{n})$ in $L^{N}(R^{n})$. In the fifth
section we give a sufficient condition which assures the strong convergence
of a sequence of solutions. In the last section, we analyze situations where
nodal solutions exist and by the way we deduce that $\mu _{2}(M,g)$ is not
attained by a classical conformal metric.

Now, we quote some facts which will be of use in the sequel of this paper.

\begin{lemma}
\label{lem1} (\cite{4}) Let $(M,g)$ be a compact Riemannian manifold of
dimension $n\geq 5$, for any $\epsilon >0$ there exists a constant $%
A(\epsilon )$ such that every $u\in H_{2}^{2}(M)$ fulfills%
\begin{equation*}
\left\Vert u\right\Vert _{N}^{2}\leq (K_{2}^{2}+\epsilon )\left\Vert \Delta
u\right\Vert _{2}^{2}+A(\epsilon )\left\Vert u\right\Vert _{2}^{2}
\end{equation*}%
with $N=\frac{2n}{n-4}$ and $K_{2}^{-2}=\pi ^{2}n(n-1)(n^{2}-4)\frac{\Gamma (%
\frac{n}{2})}{\Gamma (n)}$

where $\Gamma $ denotes the Euler function.
\end{lemma}

\begin{lemma}
\label{lem1''} ( \cite{7})Let $\left( S^{n},h\right) $ be the standard unit
sphere of $R^{n+1}$, $n\geq 5$, and let $P$ be the Paneitz-Branson operator
on $\left( S^{n},h\right) $, then%
\begin{equation*}
K_{2}^{-2}=\inf_{u\in C^{\infty }(S^{n})-\left\{ 0\right\} }\frac{%
\int_{S^{n}}uP(u)dv_{h}}{\left( \int_{S^{n}}\left\vert u\right\vert
^{N}dv_{h}\right) ^{\frac{2}{N}}}\text{.}
\end{equation*}
\end{lemma}

\begin{lemma}
\label{lem1'} ( \cite{7})Let $(M,g)$ be a smooth compact $n$-dimensional ($%
n\geq 5$) Riemannian manifold, $\alpha $ a positive real number, let $b$ be
a real valued functions defined on $M$ and $u\in H_{2}^{2}(M)$ be a weak
solution of 
\begin{equation*}
\Delta ^{2}u+\alpha \Delta u+\frac{\alpha ^{2}}{4}u=bu\text{.}
\end{equation*}%
If $b\in L^{\frac{n}{4}}(M)$, then $u\in L^{s}(M)$ for all $s\geq 1$.
\end{lemma}

\section{ First and second eigenvalues for a generalized metric}

Let $L_{+}^{N}(M)$ be the space of $L^{N}$-integrable non negative functions
which are not identically $0$. Denote by $Gr_{k}^{u}(H_{2}^{2})$ the set of
all $k$-dimensional subspaces ($k\geq 1$) of $H_{2}^{2}(M)$ which are the
span of the functions $u_{1}$, ..., $u_{k}$ if and only if $u_{1\mid
_{M-u_{1}^{-1}(0)}}$,..., $u_{k\mid _{M-u_{k}^{-1}(0)}\text{ }}$are linearly
\ independent.

\begin{definition}
A generalized metric conform to a metric $g$ is of the form $\widetilde{g}%
=ug $ with $u\in L_{+}^{N}(M)$.
\end{definition}

\begin{definition}
For any generalized metric $\widetilde{g}=u^{\frac{N-2}{2}}g$ of a
Riemannian metric $g$ we define the eigenvalue of order \ $k\geq 1$ to the
Branson-Paneitz operator $P_{g}$ by 
\begin{equation*}
\lambda _{k}(\widetilde{g})=\inf_{V\in Gr_{k}^{u}(H_{2}^{2}(M))}\sup_{v\in
V-\left\{ 0\right\} }\frac{\int_{M}vP_{g}(v)dv_{g}}{%
\int_{M}u^{N-2}v^{2}dv_{g}}\text{.}
\end{equation*}
\end{definition}

We need the following lemma which is first given in (\cite{1}) for sequences
in $H_{1}^{2}(M)$ but its proof remains inchanged and we reproduce it here
for reason of completness.

\begin{lemma}
\label{lem2} If $u\in L_{+}^{N}(M)$ and $(v_{n})$ is a sequence in $%
H_{2}^{2}(M)$ which converges weakly to $v$, then 
\begin{equation}
\int_{M}u^{N-2}\left\vert v_{m}^{2}-v^{2}\right\vert dv_{g}\rightarrow 0.
\label{2}
\end{equation}
\end{lemma}

\begin{proof}
Letting $A$ be any real positive number, we put $u_{A}=\inf (u,A)$. Then $%
(u_{A})_{A}$ is a monotone sequences which converges pointwisely almost
everywhere to $v$, so by Lebesgue monotone convergence theorem, we get%
\begin{equation*}
\int_{M}(u^{N-2}-u_{A}^{N-2})^{\frac{N}{N-2}}dv_{g}\rightarrow 0\text{.}
\end{equation*}%
On the other hand, we have%
\begin{equation*}
\int_{M}u^{N-2}\left\vert v_{m}^{2}-v^{2}\right\vert dv_{g}\leq
\int_{M}u_{A}^{N-2}\left\vert v_{m}^{2}-v^{2}\right\vert dv_{g}
\end{equation*}%
\begin{equation*}
+\int_{M}\left( u^{N-2}-u_{A}^{N-2}\right) (\left\vert v_{m}\right\vert
+\left\vert v\right\vert )^{2}dv_{g}.
\end{equation*}%
Using the H\"{o}lder inequality, we obtain%
\begin{equation*}
\int_{M}u^{N-2}\left\vert v_{m}^{2}-v^{2}\right\vert dv_{g}\leq
A^{N-2}\int_{M}\left\vert v_{m}^{2}-v^{2}\right\vert dv_{g}
\end{equation*}%
\begin{equation*}
+\left( \int_{M}\left\vert u^{N-2}-u_{A}^{N-2}\right\vert ^{\frac{N}{N-2}%
}dv_{g}\right) ^{\frac{N-2}{N}}\left( \int_{M}(\left\vert v_{m}\right\vert
+\left\vert v\right\vert )^{N}dv_{g}\right) ^{\frac{2}{N}}\text{.}
\end{equation*}%
Taking account of the boundedness and the strong convergence of $(v_{m})$ to 
$v$ in $L^{2}(M)$ we get the result.
\end{proof}

\begin{proposition}
\label{p1} Let $\widetilde{g}=u^{\frac{N-2}{2}}g$ be \ any generalized
conformal metric to a metric $g$. The equation 
\begin{equation}
P_{g}v=\lambda _{1}u^{N-2}v  \label{3}
\end{equation}%
has a solution of class $C^{4,\alpha }(M)$ ( $0<\alpha <N-2$) with the
constraint

\begin{equation*}
\int_{M}u^{N-2}v^{2}dv_{g}=1.
\end{equation*}
\end{proposition}

\begin{proof}
Let $(v_{m})$ be a minimizer sequence of $\lambda _{1}(\widetilde{g}$ ) with
the constraint $\int_{M}u^{N-2}v_{m}^{2}dv_{g}=1$. The sequence $(v_{m})$ is
bounded in $H_{2}^{2}(M)$ and by passing to a subsequences also labelled $%
(v_{m})$, there exists $v\in H_{2}^{2}(M)$ such that

(i) $v_{m}\rightarrow v$ weakly in $H_{2}^{2}(M)$

(ii) $v_{m}\rightarrow v$ strongly in $L^{2}(M)$

From (i), we obtain%
\begin{equation*}
\left\Vert v\right\Vert \leq \lim \inf \left\Vert v_{m}\right\Vert
\end{equation*}%
and by Lemma\ref{lem2} we get 
\begin{equation*}
\int_{M}u^{N-2}v^{2}dv_{g}=\lim_{m\rightarrow \infty
}\int_{M}u^{N-2}v_{m}^{2}dv_{g}=1
\end{equation*}%
and we derive that $\left\Vert v\right\Vert ^{2}=\lambda _{1}(\widetilde{g})$%
.

Consequently $v$ is a non trivial weak solution of the equation(\ref{3}).

By Lemma\ref{lem1'}, $v\in L^{s}(M)$ for any $s\geq 1$ and it follows that $%
v\in C^{4,\alpha }(M)$, with $\alpha <N-2$.
\end{proof}

\subsection{Positivity of solutions}

Now we are going to show that the equation (\ref{3}) admits a positive
solution.

\begin{proposition}
\label{p2} If the scalar curvature $S_{g}$ of the Einsteinian manifold $%
(M,g) $ is positive, the equation 
\begin{equation}
P_{g}f=\lambda _{1}u^{N-2}v  \label{4}
\end{equation}%
has a positive solution with the constraint 
\begin{equation}
\int_{M}u^{N-2}v^{2}dv_{g}=1\text{.}  \label{5}
\end{equation}
\end{proposition}

\begin{proof}
Let $v$ be a solution to the equation(\ref{4}) and let $f$ \ be the solution
of the equation 
\begin{equation*}
\Delta f+\frac{\alpha }{2}f=\left\vert \Delta v+\frac{\alpha }{2}v\right\vert
\end{equation*}%
with $\alpha >0$. Clearly $f\in C^{2,\alpha }(M)$ ($\alpha <N-2$).

If $\Delta v+\frac{\alpha }{2}v\geq 0$ ( resp. $\Delta v+\frac{\alpha }{2}%
v\leq 0$) we have $f=v$ ( resp. $f=-v$).

If it is not the case, putting $w=f\pm v$, we get 
\begin{equation}
\Delta w\pm \frac{\alpha }{2}w=\left\vert \Delta f+\frac{\alpha }{2}%
f\right\vert \pm \left( \Delta v+\frac{\alpha }{2}v\right) \geq 0  \label{6}
\end{equation}%
so $\Delta (-w)\pm \frac{\alpha }{2}\left( -w\right) \leq 0$. The maximum
principle asserts that $-w=v\pm f$ attains a maximum $M\geq 0$ then $w$ is a
constant function but this is excluded since $-\alpha $ $M\leq 0$ implies
that $M=0$. Consequently $f>\left\vert v\right\vert \geq 0$.

Let $k\geq 0$ be a real number such that $\int_{M}u^{N-2}(kf)^{2}dv_{g}=1$,
then $0<k<1$. Now letting $\widehat{f}=kf$ \ and taking account of the
equation(\ref{4}) we get%
\begin{equation*}
\int_{M}\left( \left( \Delta \widehat{f}\right) ^{2}+\alpha \left\vert
\nabla \widehat{f}\right\vert ^{2}+\overline{\alpha }\widehat{f}^{2}\right)
dv_{g}-\lambda _{1}(\widetilde{g})=
\end{equation*}%
\begin{equation*}
k^{2}\int_{M}\left( \left( \Delta f\right) ^{2}+\alpha \left\vert \nabla
f\right\vert ^{2}+\overline{\alpha }f^{2}\right) dv_{g}-\lambda _{1}(%
\widetilde{g})=
\end{equation*}%
\begin{equation*}
k^{2}\int_{M}\left( \left( \Delta f+\frac{\alpha }{2}f\right) ^{2}-\frac{%
\alpha ^{2}}{4}f^{2}+\overline{\alpha }f^{2}\right) dv_{g}-\lambda _{1}(%
\widetilde{g})=
\end{equation*}%
\begin{equation*}
k^{2}\int_{M}\left( \left( \Delta v+\frac{\alpha }{2}v\right) ^{2}-\frac{%
\alpha ^{2}}{4}f^{2}+\overline{\alpha }f^{2}\right) dv_{g}-\lambda _{1}(%
\widetilde{g})=
\end{equation*}%
\begin{equation*}
(k^{2}-1)\lambda _{1}(\widetilde{g})+\left( \overline{\alpha }-\frac{\alpha
^{2}}{4}\right) \int_{M}\left( f^{2}-v^{2}\right) dv_{g}\leq 0\text{.}
\end{equation*}%
Consequently%
\begin{equation}
\int_{M}\left( \left( \Delta \widehat{f}\right) ^{2}+\alpha \left\vert
\nabla \widehat{f}\right\vert ^{2}+\overline{\alpha }\widehat{f}^{2}\right)
dv_{g}=\lambda _{1}(\widetilde{g})\text{.}  \label{7}
\end{equation}
\end{proof}

\begin{proposition}
\label{p3} Let $u\in L_{+}^{N}(M)$, if $v\in H_{2}^{2}(M)$ is a weak
solution of the equation 
\begin{equation}
P_{g}v=\lambda _{1}(\widetilde{g})u^{N-2}v  \label{8}
\end{equation}%
with 
\begin{equation}
\int_{M}u^{N-2}v^{2}dv_{g}=1  \label{9}
\end{equation}%
then there is a weak solution $w\in H_{2}^{2}(M)$ of the equation 
\begin{equation}
P_{g}w=\lambda _{2}^{\prime }(\widetilde{g})u^{N-2}w  \label{10}
\end{equation}%
with the constraints 
\begin{equation*}
\int_{M}u^{N-2}w^{2}dv_{g}=1,\int_{M}u^{N-2}vwdv_{g}=0
\end{equation*}%
where 
\begin{equation*}
\lambda _{2}^{\prime }(\widetilde{g})=\inf_{E}\frac{\int_{M}wP_{g}wdv_{g}}{%
\int_{M}u^{N-2}w^{2}dv_{g}}
\end{equation*}%
and 
\begin{equation*}
E=\left\{ u^{\frac{N-2}{2}}w\text{ : }w\in H_{2}^{2}(M)/\ u^{\frac{N-2}{2}%
}w\ncong 0\text{, }\int_{M}u^{N-2}vwdv_{g}=0\text{ and }%
\int_{M}u^{N-2}w^{2}dv_{g}=1\text{ }\right\} \text{.}
\end{equation*}
\end{proposition}

\begin{proof}
First, we show that the set $E$ is non empty. Let $v$ , $s\in H_{2}^{2}(M)$
noncolinear such that $\int_{M}u^{N-2}v^{2}dv_{g}=1$, $%
\int_{M}u^{N-2}s^{2}dv_{g}=1$. Necessarily $u^{\frac{N-2}{2}}v\ncong 0$ and $%
u^{\frac{N-2}{2}}s\ncong 0$. \ \ Observe that $\int_{M}u^{N-2}vsdv_{g}\neq 1$%
, since if it is not the case the equality is attained in the the H\"{o}lder
inequality and this possible if and only if there a real constant $c$ such
that $v=cs$.

Putting $w=\alpha v+\beta s$ with $\alpha $, $\beta \in R$, we \ obtain%
\begin{equation*}
u^{N-2}w=\alpha u^{N-2}v+\beta u^{N-2}s
\end{equation*}

so to get 
\begin{equation*}
\int_{M}u^{N-2}vwdv_{g}=\alpha +\beta \int_{M}u^{N-2}vsdv_{g}=0
\end{equation*}%
and 
\begin{equation*}
\int_{M}u^{N-2}w^{2}dv_{g}=1
\end{equation*}%
we let 
\begin{equation*}
\beta =-\frac{\alpha }{\int_{M}u^{N-2}vsdv_{g}}
\end{equation*}%
and 
\begin{equation*}
1=\int_{M}u^{N-2}\left( \alpha v+\beta s\right) ^{2}dv_{g}
\end{equation*}%
\begin{equation*}
=\alpha ^{2}+\beta ^{2}+2\alpha \beta \int_{M}u^{N-2}vsdv_{g}\text{.}
\end{equation*}%
We obtain%
\begin{equation*}
\alpha =\pm \left( \frac{\int_{M}u^{N-2}vsdv_{g}}{1-\int_{M}u^{N-2}vsdv_{g}}%
\right) ^{\frac{1}{2}}
\end{equation*}%
and%
\begin{equation*}
\beta =\pm \frac{1}{\left( \left( 1-\int_{M}u^{N-2}vsdv_{g}\right)
\int_{M}u^{N-2}vsdv_{g}\right) ^{\frac{1}{2}}}\text{.}
\end{equation*}%
Now we will show that $w$ is a weak non trivial solution of the equation(\ref%
{10}). Let $(w_{n})$ be a minimizer sequence of $\lambda _{2}^{\prime }(%
\widetilde{g})$ such that 
\begin{equation*}
\int_{M}u^{N-2}w_{m}^{2}dv_{g}=1
\end{equation*}%
and 
\begin{equation*}
\int_{M}u^{N-2}w_{m}vdv_{g}=0\text{.}
\end{equation*}%
Then the sequences $(w_{m})$ is bounded in $H_{2}^{2}(M)$ and there is $w\in
H_{2}^{2}(M)$ a weak solution of the equation(\ref{10}). It remains to
verify that $\int_{M}u^{N-2}w^{2}=1$ and also $\int_{M}u^{N-2}wvdv_{g}=0$.
The first equality follows from Lemma\ref{lem2} the second one is true since
the function $\varphi =u^{N-2}v\in L^{\frac{N}{N-1}}(M)$.
\end{proof}

\begin{proposition}
\label{p5} Suppose that the solutions $v$ and $w$ of the equations (\ref{8})
and (\ref{10}) are as in proposition\ref{prop3}, then $\lambda _{2}(%
\widetilde{g})=\lambda _{2}^{\prime }(\widetilde{g})$.
\end{proposition}

\begin{proof}
The weak solution $w\in H_{2}^{2}(M)$ of the equation 
\begin{equation*}
P_{g}w=\lambda _{2}^{\prime }(\widetilde{g})u^{N-2}w
\end{equation*}%
is a minimizer of 
\begin{equation*}
\lambda _{2}^{\prime }(\widetilde{g})=\inf_{w\in E}\frac{%
\int_{M}wP_{g}wdv_{g}}{\int_{M}u^{N-2}w^{2}dv_{g}}
\end{equation*}%
where 
\begin{equation*}
E=\left\{ u^{\frac{N-2}{2}}w:w\in H_{2}^{2}(M)\text{ s. t. }u^{\frac{N-2}{2}%
}w\ncong 0\text{ ,}\int_{M}u^{N-2}vwdv_{g}=0\text{ and }%
\int_{M}u^{N-2}w^{2}dv_{g}=1\right\} \text{.}
\end{equation*}%
Since \ $u^{\frac{N-2}{2}}v$ and $u^{\frac{N-2}{2}}w$ are linearly
independent it follows that $V_{o}=$span($v,w)\in Gr_{2}^{u}(H_{2}^{2}(M))$.

Putting 
\begin{equation*}
f=\lambda v+\mu w\text{ \ with }(\lambda ,\mu )\in R^{2}-\left\{ \left(
0,0\right) \right\}
\end{equation*}
we evaluate 
\begin{equation*}
s=\frac{\int_{M}fP_{g}fdv_{g}}{\int_{M}u^{N-2}f^{2}dv_{g}}
\end{equation*}%
on the plane $V_{o}$.

We obtain%
\begin{equation*}
s=\frac{\lambda ^{2}\int_{M}vP_{g}\left( v\right) dv_{g}+\mu
^{2}\int_{M}wP_{g}\left( w\right) dv_{g}}{\lambda ^{2}+\mu ^{2}}
\end{equation*}%
\begin{equation*}
=\frac{\lambda ^{2}}{\lambda ^{2}+\mu ^{2}}\lambda _{1}(\widetilde{g})+\frac{%
\mu ^{2}}{\lambda ^{2}+\mu ^{2}}\lambda _{2}^{\prime }(\widetilde{g})
\end{equation*}%
\begin{equation*}
=\lambda _{1}(\widetilde{g})\cos ^{2}\theta +\lambda _{2}^{\prime }(%
\widetilde{g})\sin ^{2}\theta
\end{equation*}%
with $\theta \in \mathbb{R}$.

On the other hand, we have%
\begin{equation*}
\frac{ds}{d\theta }=(\lambda _{2}^{\prime }(\widetilde{g})-\lambda _{1}(%
\widetilde{g}))\sin 2\theta
\end{equation*}%
and noting that

\begin{equation*}
\lambda _{1}(\widetilde{g})\leq \lambda _{2}^{\prime }(\widetilde{g})
\end{equation*}%
we get easily 
\begin{equation*}
\min s(\theta )=\lambda _{1}(\widetilde{g})\text{ and }\max s(\theta
)=\lambda _{2}^{\prime }(\widetilde{g})\text{.}
\end{equation*}%
Consequently%
\begin{equation*}
\lambda _{2}^{\prime }(\widetilde{g})=\sup_{w\in V_{o}}\frac{%
\int_{M}wP_{g}(w)dv_{g}}{\int_{M}u^{N-2}w^{2}dv_{g}}\text{.}
\end{equation*}%
On the other hand the infimum of $\sup_{w\in V-\left\{ 0\right\} }\frac{%
\int_{M}wP_{g}(w)dv_{g}}{\int_{M}u^{N-2}w^{2}dv_{g}}$ on all the subspaces
of $Gr_{2}^{u}(H_{2}^{2}(M))$ is attained by $V_{o}=span(v,w)$. \ 

Hence%
\begin{equation*}
\lambda _{2}^{\prime }(\widetilde{g})=\lambda _{2}(\widetilde{g})\text{.}
\end{equation*}%
{}
\end{proof}

\begin{proposition}
\label{p4} If $u\in C^{\infty }(M)$ with $u\geq 0$ not identically $0$. Then
any weak solution of the equation 
\begin{equation}
P_{g}v=\mu u^{N-2}v  \label{11}
\end{equation}%
is of class $C^{\infty }(M)$, $\mu \in \mathbb{R}$.
\end{proposition}

\begin{proof}
Let $u\in C^{\infty }(M)$, $u\geq 0$ and not identically $0$ and $v$ a weak
solution of the equation(\ref{11}). We have%
\begin{equation*}
(\Delta +a)(\Delta +b)v=\mu u^{N-2}v
\end{equation*}%
with $a=\frac{\alpha -\sqrt{\alpha ^{2}-4\overline{\alpha }}}{2}$ and $b=%
\frac{\alpha +\sqrt{\alpha ^{2}-4\overline{\alpha }}}{2}$.

Putting 
\begin{equation*}
z=(\Delta +b)v
\end{equation*}%
we get%
\begin{equation*}
(\Delta +a)z=\mu u^{N-2}v
\end{equation*}%
and since $v\in H_{2}^{2}(M)$, $u^{N-2}v\in H_{2}^{2}(M)$ so $z\in
H_{4}^{2}(M)$. Recurrently, for any $k\geq 2$ we obtain $v\in H_{k}^{2}$.
Now, classical regularity theorem allows us to conclude that $v\in C^{\infty
}(M)$.
\end{proof}

\section{A Sobolev inequality related to $\protect\mu _{2}(M,g)$}

The Sobolev inequality given by Lemma\ref{lem1} which allows to avoid
concentration phenomena for the minimizing sequence of the first
Paneitz-Branson invaiant $\mu _{1}(M,g)$ is not sufficient in the case of
the second Paneitz-Branson invariant $\mu _{2}(M,g)$, we propose the
following Sobolev type inequality.

\begin{proposition}
\label{prop3} Let $(M,g)$ be a Riemannian manifold of dimension $n\geq 5$.
For any $\epsilon >0$ there is a constant $A(\epsilon )$ such that, for any $%
u\in L_{+}^{N}(M)$ and any $v\in H_{2}^{2}(M)$, we have 
\begin{equation*}
\int_{M}u^{N-2}v^{2}dv_{g}\leq \left( 2^{-\frac{4}{n}}(K_{2}^{2}+\epsilon
)\int_{M}(\Delta v)^{2}dv_{g}+A(\epsilon )\int_{M}v^{2}dv_{g}\right) \left(
\int_{M}u^{N}dv_{g}\right) ^{\frac{2}{N}}
\end{equation*}
\end{proposition}

\begin{proof}
For any $\epsilon >0$, put 
\begin{equation*}
B(\epsilon )=A(\epsilon )K_{2}^{-2}(1+\epsilon )^{-1}
\end{equation*}%
and let 
\begin{equation*}
G(u,v)=\frac{\int_{M}(\Delta v)^{2}dv_{g}+B(\epsilon )\int_{M}v^{2}dv_{g}}{%
\int_{M}u^{N-2}v^{2}dv_{g}}\left( \int_{M}u^{N}dv_{g}\right)
\end{equation*}%
where $u\in L_{+}^{N}(M)$ and $v\in H_{2}^{2}(M)-\left\{ 0\right\} $ such
that $\int_{M}u^{N-2}v^{2}dv_{g}\neq 0$.

Obviously $G(u,v)$ is continuous on $L_{+}^{N}(M)\times H_{2}^{2}(M)-\left\{
0\right\} $. So $I(u,V)=\sup_{v\in V-\left\{ 0\right\} }G(u,v)$ depends
continuously on $u\in L_{+}^{N}(M)$ and $V\in Gr_{2}^{u}(H_{2}^{2}(M))$. We
must show that 
\begin{equation*}
I(u,V)\geq 2^{\frac{4}{n}}K_{2}^{-2}(1+\epsilon )^{-1}
\end{equation*}%
for all $u\in C^{\infty }(M)$, $u>0$ and $V\in Gr_{2}^{u}(C^{\infty }(M))$.

Without lost of generality, we suppose that $\int_{M}u^{N-2}v^{2}dv_{g}=1$.
On the other hand the operator 
\begin{equation}
v\rightarrow Q(v)=u^{\frac{2-N}{2}}\Delta ^{2}(u^{\frac{2-N}{2}%
}v)+B(\epsilon )u^{2-N}v  \label{Q}
\end{equation}%
is a fourth order elliptic and self adjoint with respect to the inner
product in $L^{2}(M)$. $Q$ has a discrete spectrum $\ \lambda _{1}\leq
\lambda _{2}\leq ...$

The corresponding eigenfunctions $\varphi _{1},\varphi _{2},...$ are smooth
functions on $M$. Letting $v_{i}=u^{\frac{2-N}{2}}\varphi _{i}$ , we get 
\begin{equation*}
\int_{M}(\Delta v_{i})^{2}dv_{g}+B(\epsilon )\int_{M}v_{i}^{2}dv_{g}=\lambda
_{i}\int_{M}u^{N-2}v_{i}^{2}dv_{g}
\end{equation*}%
with 
\begin{equation*}
\int_{M}u^{N-2}v_{i}v_{j}dv_{g}=0\text{.}
\end{equation*}

Let $\widetilde{P}_{g}$ be the operator defined on $C^{\infty }(M)$ by $%
\widetilde{P}_{g}u=\Delta _{g}^{2}u+B(\epsilon )u$ and let $\Omega _{1}$ and 
$\Omega _{2}$ be two non empty open disjoint sets in $M$ and let $v_{1}$ and 
$v_{2}$ be two non trivial solutions to the equation 
\begin{equation}
\widetilde{P}_{g}v_{i}=\lambda _{2}u^{N-2}v_{i}  \label{E}
\end{equation}%
$i=1,2$ with supports included respedtively in $\overline{\Omega }_{1}$ and $%
\overline{\Omega }_{2}$, the closer sets of $\Omega _{1}$ and $\Omega _{2}$
r and where $\lambda _{2}$ is the second eigenvalue of the operator $Q$
defined above. By multiplying if necessary $v_{1}$and $v_{2}$ by constants,
we assume that $\int_{M}u^{N-2}v_{1}^{2}dv_{g}=%
\int_{M}u^{N-2}v_{2}^{2}dv_{g}=1$.

Using the H\"{o}lder inequality and the Sobolev one given in Lemma\ref{lem1}%
, we get%
\begin{equation*}
2=\int_{M}u^{N-2}v_{1}^{2}dv_{g}+\int_{M}u^{N-2}v_{2}^{2}dv_{g}
\end{equation*}%
\begin{equation*}
\leq \left( \int_{\Omega _{1}}u^{N}dv_{g}\right) ^{1-\frac{2}{N}}\left(
\int_{M}\left\vert v_{1}\right\vert ^{N}dv_{g}\right) ^{\frac{2}{N}}+\left(
\int_{\Omega _{2}}u^{N}dv_{g}\right) ^{1-\frac{2}{N}}\left(
\int_{M}\left\vert v_{2}\right\vert ^{N}dv_{g}\right) ^{\frac{2}{N}}
\end{equation*}

\begin{equation*}
\leq \left( \int_{\Omega _{1}}u^{N}dv_{g}\right) ^{1-\frac{2}{N}%
}K_{2}^{2}(1+\epsilon )\left( \int_{M}(\Delta v_{1})^{2}dv_{g}+B(\epsilon
)\int_{M}v_{1}^{2}dv_{g}\right)
\end{equation*}%
\begin{equation*}
+\left( \int_{\Omega _{2}}u^{N}dv_{g}\right) ^{1-\frac{2}{N}%
}K_{2}^{2}(1+\epsilon )\left( \int_{M}(\Delta v_{2})^{2}dv_{g}+B(\epsilon
)\int_{M}v_{2}^{2}dv_{g}\right) .
\end{equation*}%
And since $v_{1}$ and $v_{2}$ are solutions to the equation (\ref{E}), we
obtain%
\begin{equation*}
2\leq K_{2}^{2}(1+\epsilon )\lambda _{2}\left( \left( \int_{\Omega
_{1}}u^{N}dv_{g}\right) ^{1-\frac{2}{N}}+\left( \int_{\Omega
_{2}}u^{N}dv_{g}\right) ^{1-\frac{2}{N}}\right) \text{.}
\end{equation*}%
Using H\"{o}lder inequality, we get 
\begin{equation*}
\left( \int_{\Omega _{1}}u^{N}dv_{g}\right) ^{1-\frac{2}{N}}+\left(
\int_{\Omega _{2}}u^{N}dv_{g}\right) ^{1-\frac{2}{N}}\leq 2^{\frac{2}{N}%
}\left( \int_{\Omega _{1}}u^{N}dv_{g}+\int_{\Omega _{2}}u^{N}dv_{g}\right)
\end{equation*}%
so 
\begin{equation*}
\lambda _{2}\geq 2^{\frac{4}{n}}(1+\epsilon )^{-1}K_{2}^{-2}\text{.}
\end{equation*}%
Letting $V=span(v_{1},v_{2})$, we obtain for any $(\alpha ,\beta )\in
R^{2}-(0,0)$,

\begin{equation*}
G(u,\alpha v_{1}+\beta v_{2})=\frac{\int_{M}\left[ \left( \Delta (\alpha
v_{1}+\beta v_{2})\right) ^{2}+B(\epsilon )(\alpha v_{1}+\beta v_{2})^{2}%
\right] dv_{g}}{\int_{M}u^{n-2}(\alpha v_{1}+\beta v_{2})^{2}dv_{g}}
\end{equation*}%
\begin{equation*}
=\frac{\alpha ^{2}\int_{M}\left( (\Delta v_{1})^{2}+B(\epsilon
)v_{1}^{2}\right) dv_{g}+\alpha ^{2}\int_{M}\left( (\Delta
v_{2})^{2}+B(\epsilon )v_{2}^{2}\right) dv_{g}}{\alpha
^{2}\int_{M}u^{N-2}v_{1}^{2}dv_{g}+\beta ^{2}\int_{M}u^{N-2}v_{2}^{2}dv_{g}}
\end{equation*}%
\begin{equation*}
=\lambda _{2}\text{.}
\end{equation*}%
Then 
\begin{equation*}
I(u,V)=\sup_{(\alpha ,\beta )\in R^{2}-(0,0)}G(u,\alpha v_{1}+\beta
v_{2})=\lambda _{2}
\end{equation*}%
and the proof of the proposition is achieved.
\end{proof}

In the particular case of the standard unit $\left( S^{n},h\right) $ sphere
of $R^{n+1}$, we obtain

\begin{proposition}
\label{prop3'}Let $\left( S^{n},h\right) $ be the unit sphere of $R^{n+1}$, $%
n\geq 5$, and let $P$ be the Paneitz-Branson operator on $\left(
S^{n},h\right) $ . For any $u\in L_{+}^{N}(S^{n})$ and any $v\in
H_{2}^{2}(S^{n})$, we have 
\begin{equation*}
\int_{S^{n}}u^{N-2}v^{2}dv_{h}\leq 2^{-\frac{4}{n}}K_{2}^{2}%
\int_{S^{n}}vP(v)dv_{h}\left( \int_{S^{n}}u^{N}dv_{h}\right) ^{\frac{2}{N}}%
\text{.}
\end{equation*}
\end{proposition}

\begin{proof}
The proof is similar to that of the propostion\ref{prop3}, by using the
Sobolev inequality given by Lemma\ref{lem1''} instead of that given by Lemma%
\ref{lem1}.
\end{proof}

As corollary of proposition\ref{prop3'}, we get the following Sobolev
inequality on the Euclidean space $R^{n}.$

\begin{corollary}
\label{cor1} Let $C_{c}^{\infty }(R^{n})$ be the space of functions of
classe $C^{\infty }$ and of compact supports on $R^{n}$. For any $u\in
L_{+}^{N}(R^{n})$ and any $v\in H_{2}^{2}(R^{n})$, we have 
\begin{equation*}
\int_{R^{n}}u^{N-2}v^{2}dx\leq 2^{-\frac{4}{n}}K_{2}^{2}\int_{R^{n}}(\Delta
v)^{2}dx\left( \int_{M}u^{N}dx\right) ^{\frac{2}{N}}
\end{equation*}%
where $dx$ denotes the Euclidean measure on $R^{n}$.
\end{corollary}

\begin{proof}
Since $R^{n}$ is conformal to $S^{n}-\left\{ p\right\} $, where $p$ is any
point of $S^{n}$ and the Paneitz-Branson is a conformal invariant the
corollory\ref{cor1} follows from proposition\ref{prop3'}.
\end{proof}

\begin{proposition}
\bigskip If $\mu _{1}(M,g)K_{2}^{2}<1$, then%
\begin{equation*}
\mu (M,g)=\mu _{1}(M,g)
\end{equation*}
\end{proposition}

\begin{proof}
\begin{equation*}
\mu _{1}(M,g)=\inf_{\widetilde{g}\in \left[ g\right] }\lambda _{1}(%
\widetilde{g})\left( vol(M)\right) ^{\frac{4}{n}}
\end{equation*}%
\begin{equation*}
=\inf_{\substack{ u\in C^{\infty }(M)  \\ u>0}}\inf_{v\in C^{\infty
}(M)-\left\{ 0\right\} }\frac{\int_{M}vP_{g}(v)dv_{g}}{%
\int_{M}u^{N-2}v^{2}dv_{g}}\left( \int_{M}u^{N}dv_{g}\right) ^{\frac{4}{n}}
\end{equation*}%
\begin{equation*}
\leq \inf_{v\in C^{\infty }(M)-\left\{ 0\right\} }\frac{%
\int_{M}vP_{g}(v)dv_{g}}{\left( \int_{M}\left\vert v\right\vert
^{N}dv_{g}\right) ^{\frac{2}{N}}}=\mu (M,g)\text{.}
\end{equation*}%
The inequality in the other sense requires a variational method. Let $%
g_{m}=u_{m}^{\frac{N-2}{2}}g$ \ with $u_{m}\in L_{+}^{N}(M)$, a minimizer
sequence of $\mu _{1}(M,g)$ i.e.%
\begin{equation*}
\mu _{1}(M,g)=\lim_{m\rightarrow \infty }\lambda _{1}(g_{m})\left(
vol(M,g)\right) ^{\frac{4}{n}}\text{.}
\end{equation*}

Considering the Yamabe functional%
\begin{equation*}
Y(u,v)=\frac{\int_{M}vP_{g}(v)dv_{g}}{\int_{M}u^{N-2}v^{2}dv_{g}}\left(
\int_{M}u^{N}dv_{g}\right) ^{\frac{4}{n}}
\end{equation*}%
with $v\in H_{2}^{2}(M)-\{0\}$ and $u\in L_{+}^{N}(M)$, we write, for any $%
\lambda \in \mathbb{R}^{\ast }$%
\begin{equation*}
Y(\lambda u,v)=\frac{\int_{M}vP_{g}(v)dv_{g}}{\lambda
^{N-2}\int_{M}u^{N-2}v^{2}dv_{g}}\lambda ^{\frac{2N}{n}}\left(
\int_{M}u^{N}dv_{g}\right) ^{\frac{4}{n}}
\end{equation*}%
\begin{equation*}
=Y(u,v)\text{.}
\end{equation*}%
So, we can choose the sequence $(u_{m})$ such that $%
\int_{M}u_{m}^{N}dv_{g}=1 $ and there is a subsequence of $(u_{m})$ still
labelled by $(u_{m})$ converging weakly to $u\geq 0$ in $L^{N}(M)$.

On the other hand, by Proposition\ref{p3} for any $u_{m}\in L_{+}^{N}(M)$
there is $v_{m}\in H_{2}^{2}(M)$ solutions of the equation%
\begin{equation*}
P_{g}(v_{m})=\lambda _{1,m}u_{m}^{N-2}v_{m}
\end{equation*}%
with the constraint%
\begin{equation*}
\int_{M}u_{m}^{N-2}v_{m}^{2}dv_{g}=1\text{.}
\end{equation*}%
Obviously $(v_{m})$ is bounded in $H_{2}^{2}(M)$ so there is $v\in
H_{2}^{2}(M)$ such that $v_{m}\rightarrow v$ \ weakly in $H_{2}^{2}(M)$, $%
v_{m}\rightarrow v$ a.e. in $M$.

Since $\lim_{m\rightarrow \infty }\lambda _{1,m}=\mu _{1}(M,g)$, $v$ is a
weak solution of the equation 
\begin{equation*}
P_{g}(v)=\mu _{1}(M,g)u^{N-2}v\text{ .}
\end{equation*}

Now, we are going to show that $v$ satisfies the condition%
\begin{equation*}
\int_{M}u^{N-2}v^{2}dv_{g}=1\text{.}
\end{equation*}%
It is obvious that 
\begin{equation*}
\int_{M}u^{N-2}v^{2}dv_{g}\leq 1
\end{equation*}%
So we have to show the inequality in the other sense, to do so, we consider%
\begin{equation*}
\int_{M}u^{N-2}v^{2}dv_{g}=\int_{M}u_{m}^{N-2}v_{m}^{2}dv_{g}-%
\int_{M}(u_{m}^{N-2}v_{m}^{2}-u^{N-2}v^{2})dv_{g}
\end{equation*}%
\begin{equation*}
=1-\int_{M}(u_{m}^{N-2}v_{m}^{2}-u^{N-2}v^{2})dv_{g}\text{.}
\end{equation*}%
Now, since 
\begin{equation*}
\left\vert u_{m}^{N-2}v_{m}^{2}-u_{m}^{N-2}(v_{m}-v)^{2}\right\vert \leq
Cu_{m}^{N-2}\left\vert v_{m}+v\right\vert \left\vert v\right\vert
\end{equation*}%
where $C$ is a postive constant

we get%
\begin{equation*}
u_{m}^{N-2}v_{m}^{2}-u_{m}^{N-2}(v_{m}-v_{2})^{2}\rightarrow u^{N-2}v^{2}%
\text{ in }L^{1}(M)
\end{equation*}%
and%
\begin{equation*}
\int_{M}\left( u_{m}^{N-2}v_{m}^{2}-u^{N-2}v^{2}\right) dv_{g}\rightarrow
\int_{M}u_{m}^{N-2}(v_{m}-v_{2})^{2}dv_{g}\text{.}
\end{equation*}%
so%
\begin{equation}
\int_{M}u^{N-2}v^{2}dv_{g}=1-\int_{M}u_{m}^{N-2}(v_{m}-v)^{2}dv_{g}+o(1)%
\text{.}  \label{16}
\end{equation}%
where $o(1)$ is a sequence converging to $0$ as $m\rightarrow +\infty $.
Using simultaneously the H\"{o}lder inequality and the Sobolev inequality
given by Lemma\ref{lem1}, we get%
\begin{equation*}
\int_{M}u_{m}^{N-2}(v_{m}-v)^{2}dv_{g}\leq \left( K_{2}^{2}+\epsilon \right)
\left\Vert \Delta \left( v_{m}-v\right) \right\Vert _{2}^{2}+A(\epsilon
)\left\Vert v_{m}-v\right\Vert _{2}^{2}
\end{equation*}%
\begin{equation*}
\leq \left( K_{2}^{2}+\epsilon \right) \left( \int_{M}\left(
v_{m}P_{g}(v_{m})-vP_{g}(v)\right) dv_{g}\right) +o(1)
\end{equation*}%
\begin{equation*}
\leq \left( K_{2}^{2}+\epsilon \right) \mu
_{1}(M,g)\int_{M}(u_{m}^{N-2}v_{m}^{2}-u^{N-2}v^{2})dv_{g}+o(1)
\end{equation*}%
\begin{equation*}
\leq \left( K_{2}^{2}+\epsilon \right) \mu _{1}(M,g)\left(
1-\int_{M}u^{N-2}v^{2}dv_{g}\right) +o(1)\text{. }
\end{equation*}%
Taking account of (\ref{16}), we have%
\begin{equation*}
\int_{M}u^{N-2}v^{2}dv_{g}\geq 1-\left( K_{2}^{2}+\epsilon \right) \mu
_{1}(M,g)\left( 1-\int_{M}u^{N-2}v^{2}dv_{g}\right) 
\end{equation*}%
\begin{equation*}
+o(1)\text{.}
\end{equation*}%
Then%
\begin{equation*}
\left( 1-\left( K_{2}^{2}+\epsilon \right) \mu _{1}(M,g)\right)
\int_{M}u^{N-2}v^{2}dv_{g}\geq 1-2^{-\frac{4}{n}}\left( K_{2}^{2}+\epsilon
\right) \mu _{1}(M,g)
\end{equation*}%
\begin{equation*}
+o(1)\text{.}
\end{equation*}%
So if 
\begin{equation*}
\mu _{1}(M,g)K_{2}^{2}<1
\end{equation*}%
we get%
\begin{equation*}
\int_{M}u^{N-2}v^{2}dv_{g}\geq 1\text{.}
\end{equation*}%
Consequently 
\begin{equation*}
\mu _{1}(M,g)>0\text{.}
\end{equation*}%
Let $\overline{u}=a\left\vert v\right\vert $ with $a>0$ and%
\begin{equation*}
\int_{M}\overline{u}^{N}dv_{g}=a^{N}\int_{M}v^{N}dv_{g}=1
\end{equation*}%
then%
\begin{equation*}
\mu _{1}(M,g)\leq \frac{\int_{M}vP_{g}(v)dv_{g}}{\int_{M}\overline{u}%
^{N-2}v^{2}dv_{g}}
\end{equation*}%
\begin{equation*}
\leq \frac{\mu _{1}(M,g)\int_{M}u^{N-2}v^{2}dv_{g}}{a^{N-2}%
\int_{M}v^{N}dv_{g}}
\end{equation*}%
\begin{equation*}
\leq a^{2}\mu _{1}(M,g)\int_{M}u^{N-2}v^{2}dv_{g}
\end{equation*}%
\begin{equation*}
\leq \mu _{1}(M,g)\int_{M}u^{N-2}\overline{u}^{2}dv_{g}
\end{equation*}%
The H\"{o}lder inequality implies that%
\begin{equation*}
\mu _{1}(M,g)\leq \mu _{1}(M,g)\left( \int_{M}u^{N}dv_{g}\right) ^{1-\frac{2%
}{N}}\left( \int_{M}\overline{u}^{N}dv_{g}\right) ^{\frac{2}{N}}
\end{equation*}%
\begin{equation*}
\leq \mu _{1}(M,g)\text{.}
\end{equation*}%
So the equality is attained in the H\"{o}lder inequality and this is
possible only if 
\begin{equation*}
\overline{u}=cu
\end{equation*}%
with $c>0$ is a constant which implies that 
\begin{equation*}
c=1
\end{equation*}%
and 
\begin{equation*}
u=\overline{u}=a\left\vert v\right\vert \text{.}
\end{equation*}%
Also%
\begin{equation*}
a^{N-2}\int_{M}v^{N}dv_{g}=\int_{M}u^{N-2}v^{2}dv_{g}=1
\end{equation*}%
and%
\begin{equation*}
\frac{a^{N}\int_{M}v^{N}dv_{g}}{a^{2}}=1
\end{equation*}%
Finally since $a>0$, we get 
\begin{equation*}
a=1
\end{equation*}%
hence%
\begin{equation*}
u=\left\vert v\right\vert \text{.}
\end{equation*}%
That means that $v$ is a weak solution in $H_{2}^{2}(M)$ to the equation 
\begin{equation*}
P_{g}(v)=\mu _{1}(M,g)\left\vert v\right\vert ^{N-2}v\text{.}
\end{equation*}%
The condition \ $\int_{M}v^{N}dv_{g}=1$ implies that $v$ is non trivial.
Consequently%
\begin{equation*}
\mu _{1}(M,g)=\frac{\int_{M}P_{g}(v)dv_{g}}{\int_{M}\left\vert v\right\vert
^{N}dv_{g}}\geq \mu (M,g)\text{.}
\end{equation*}
\end{proof}

\section{Existence of a minimizer to $\protect\mu _{2}(M,g)$}

\begin{proposition}
\label{pro1} If $\ \mu _{2}(M,g)K_{2}^{2}2^{-\frac{4}{n}}<1$, then $\mu
_{2}(M,g)$ is attained by a generalized metric $u^{N-2}g$, u$\in L_{+}^{N}(M)
$.
\end{proposition}

\begin{proof}
Let $g_{m}=u_{m}^{N-2}g$ with $u_{m}\in C^{\infty }(M)$ and $u_{m}>0$ be a
minimizing sequence of $\mu _{2}(M,g)$. Since we can assume that 
\begin{equation*}
\int_{M}u_{m}^{N}dv_{g}=1
\end{equation*}%
we have 
\begin{equation*}
\lim_{n}\lambda _{2,m}=\mu _{2}(M,g)
\end{equation*}%
By Proposition\ref{p3}, there are $v_{m},w_{m}\in H_{2}^{2}(M)$ such that 
\begin{equation}
P(v_{m})=\lambda _{1,m}u_{m}^{N-2}v_{m}  \label{17}
\end{equation}%
and 
\begin{equation}
P(w_{m})=\lambda _{2,m}u_{m}^{N-2}w_{m}  \label{18}
\end{equation}%
with the normalized conditions 
\begin{equation}
\int_{M}u_{m}^{N-2}v_{m}^{2}dv_{g}=\int_{M}u_{m}^{N-2}v_{m}^{2}dv_{g}=1\text{%
,}\ \ \int_{M}u_{m}^{N-2}v_{m}w_{m}dv_{g}=0\text{.}  \label{19}
\end{equation}

First, we have for any integer $m\geq 1$,%
\begin{equation*}
\lambda _{1,m}<\lambda _{2,m}\text{.}
\end{equation*}%
Since, if $\lambda _{1,m}=\lambda _{2,m}$; $w_{m}$ is a minimizer of $%
\lambda _{1,m}$. On other hand taking account of the coerciveness of the
Paneitz operator $P$ \ and applying the Lax-Milgram theorem, we get easily
that the first eigenvalue $\lambda _{1,m}$ of $P_{g}$ is simple, so $%
w_{m}=\alpha v_{m}$ with a real $\alpha \neq 0.$ Thus by (\ref{19}), we get
that%
\begin{equation*}
\int_{M}u_{m}^{N-2}v_{m}^{2}dv_{g}=0
\end{equation*}%
a contradiction. The sequences $(v_{m})_{m}$ and $(w_{m})_{m}$ are bounded
in $H_{2}^{2}(M)$ so there are functions $v,w\in H_{2}^{2}(M)$ \ and
subsequences still denoted by $(v_{m})_{m}$ and $(w_{m})_{m}$ converging
weakly to $v$ and $w$, respectively, in $H_{2}^{2}(M)$. This latter facts
and the weak convergence of $(u_{m})$ to $u$ in $L^{N}(M)$ allow us to write
in the weak sense that 
\begin{equation}
P_{g}(v)=\nu u^{N-2}v  \label{20}
\end{equation}%
and w%
\begin{equation}
P_{g}(w)=\mu _{2}(M,g)u^{N-2}w  \label{21}
\end{equation}%
where $\nu =\lim_{m}\lambda _{1,m}\leq \mu _{2}(M,g)$.

Now, we are going to show that $v$ ,$w$ fulfill respectively the conditions $%
\int_{M}u^{N-2}v^{2}dv_{g}=\int_{M}u^{N-2}w^{2}dv_{g}=1$ and by the way $v$, 
$w$ are not identically null. To do so, we borrow ideas and notations from(%
\cite{1}).\ Set%
\begin{equation*}
S_{m}=\left\{ \lambda _{m}v_{m}+\mu _{m}w_{m}\text{:\ }\left( \lambda
_{m},\mu _{m}\right) \in R^{2}\text{, }\lambda _{m}^{2}+\mu _{m}^{2}=1\text{%
, }\lambda _{m}\mu _{m}>\alpha >0\text{\ }\right\}
\end{equation*}%
\begin{equation*}
S=\left\{ \lambda v+\mu w\text{:\ }\left( \lambda ,\mu \right) \in R^{2}%
\text{, }\lambda ^{2}+\mu ^{2}=1\right\}
\end{equation*}%
and let $\overline{w}_{m}=\lambda _{m}v_{m}+\mu _{m}w_{m}$ , $\overline{w}%
=\lambda v+\mu w$ where up to a subsequence $\left( \lambda _{m},\mu
_{m}\right) \rightarrow \left( \lambda ,\mu \right) $.

Obviously, we have%
\begin{equation}
\int_{M}u^{N-2}\overline{w}^{2}dv_{g}\leq \lim_{m}\inf \int_{M}u_{m}^{N-2}%
\overline{w}_{m}^{2}dv_{g}=1\text{.}  \label{21'}
\end{equation}

For the inequality in the other sense, we have

\begin{equation*}
\int_{M}u^{N-2}\overline{w}^{2}dv_{g}=\int_{M}u_{m}^{N-2}\overline{w}%
_{m}^{2}dv_{g}-\int_{M}(u_{m}^{N-2}\overline{w}_{m}^{2}-u^{N-2}\overline{w}%
^{2})dv_{g}
\end{equation*}%
\begin{equation*}
=1-\int_{M}(u_{m}^{N-2}v_{m}^{2}-u^{N-2}v^{2})dv_{g}.
\end{equation*}%
Now, since 
\begin{equation}
\left\vert u_{m}^{N-2}\overline{w}_{m}^{2}-u_{m}^{N-2}(\overline{w}_{m}-%
\overline{w})^{2}\right\vert \leq Cu_{m}^{N-2}\left\vert \overline{w}_{m}+%
\overline{w}\right\vert \left\vert \overline{w}\right\vert   \label{22}
\end{equation}%
where $C>0$ is some constant

we get%
\begin{equation*}
\left\vert u_{m}^{N-2}\overline{w}_{m}^{2}-u_{m}^{N-2}(\overline{w}_{m}-%
\overline{w})^{2}\right\vert \rightarrow u^{N-2}\overline{w}^{2}\text{ in }%
L^{1}(M)
\end{equation*}%
and%
\begin{equation}
\int_{M}\left( u_{m}^{N-2}\overline{w}_{m}^{2}-u^{N-2}\overline{w}%
^{2}\right) dv_{g}\rightarrow \int_{M}u_{m}^{N-2}(\overline{w}_{m}-\overline{%
w})^{2}dv_{g}\text{.}  \label{23}
\end{equation}%
Thus%
\begin{equation*}
\int_{M}u^{N-2}\overline{w}^{2}dv_{g}=1-\int_{M}u_{m}^{N-2}(\overline{w}_{m}-%
\overline{w})^{2}dv_{g}+o(1)
\end{equation*}%
where $o(1)$ is a sequence which goes to $0$ as $m\rightarrow +\infty $.

Using the Sobolev inequality given by proposition(\ref{prop3}), and taking
account of%
\begin{equation*}
\left\Vert u_{m}\right\Vert _{N}=1
\end{equation*}%
we get 
\begin{equation*}
\int_{M}u_{m}^{N-2}(\overline{w}_{m}-\overline{w})^{2}dv_{g}\leq 2^{-\frac{4%
}{n}}\left( K_{2}^{2}+\epsilon \right) \left\Vert \Delta (\overline{w}_{m}-%
\overline{w})\right\Vert _{2}^{2}+A(\epsilon )\left\Vert \overline{w}_{m}-%
\overline{w}\right\Vert _{2}^{2}\text{.}
\end{equation*}%
Now by the Brezis-Lieb lemma(\cite{3}) and the fact that $\left\Vert 
\overline{w}_{m}-\overline{w}\right\Vert _{2}\rightarrow 0$ as $m\rightarrow
+\infty $, we obtain%
\begin{equation}
\int_{M}u_{m}^{N-2}(\overline{w}_{m}-\overline{w})^{2}dv_{g}\leq 2^{-\frac{4%
}{n}}\left( K_{2}^{2}+\epsilon \right) \left( \left\Vert \Delta \overline{w}%
_{m}\right\Vert _{2}^{2}-\left\Vert \Delta \overline{w}\right\Vert
_{2}^{2}\right) +o(1)\text{.}  \label{23'}
\end{equation}

By the fact that 
\begin{equation*}
\overline{w}_{m}-\overline{w}\rightarrow 0\text{ in }H_{q}^{2}(M)\text{, }%
q=0,1\text{ as }m\rightarrow +\infty 
\end{equation*}%
we have%
\begin{equation*}
\left\Vert \Delta \overline{w}_{m}\right\Vert _{2}^{2}-\left\Vert \Delta 
\overline{w}\right\Vert _{2}^{2}=\int_{M}\left( \overline{w}_{m}P_{g}(%
\overline{w}_{m})-\overline{w}P(\overline{w})\right) dv_{g}+o(1)\text{.}
\end{equation*}%
\begin{equation*}
=\lambda _{1,m}\lambda _{m}^{2}\int_{M}u_{m}^{N-2}v_{m}^{2}+\lambda
_{2,m}\mu _{m}^{2}\int_{M}u_{m}^{N-2}w_{m}^{2}-\nu (M,g)\lambda
^{2}\int_{M}u^{N-2}v^{2}-\mu _{2}\mu ^{2}\int_{M}u^{N-2}w^{2}+o(1)
\end{equation*}%
Taking into account of (\ref{21'}), we get%
\begin{equation*}
\left\Vert \Delta \overline{w}_{m}\right\Vert _{2}^{2}-\left\Vert \Delta 
\overline{w}\right\Vert _{2}^{2}\leq \lambda _{2,m}\left( \int_{M}u_{m}^{N-2}%
\overline{w}_{m}^{2}-\int_{M}u^{N-2}\overline{w}^{2}dv_{g}\right) 
\end{equation*}%
\begin{equation*}
+\left( \lambda _{1,m}\lambda _{m}^{2}-\nu (M,g)\lambda ^{2}\right)
\int_{M}u^{N-2}v^{2}+\left( \lambda _{2,m}\mu _{m}^{2}-\mu _{2}\mu
^{2}\right) \int_{M}u^{N-2}w^{2}+o(1)
\end{equation*}%
\begin{equation*}
\leq \mu _{2}(M,g)\int_{M}(u_{m}^{N-2}\overline{w}_{m}^{2}-u^{N-2}\overline{w%
}^{2})dv_{g}+o(1)\text{.}
\end{equation*}%
Consequently%
\begin{equation*}
\int_{M}u^{N-2}\overline{w}^{2}dv_{g}\geq 1-2^{-\frac{4}{n}}\left(
K_{2}^{2}+\epsilon \right) \mu _{2}(M,g)\int_{M}(u_{m}^{N-2}\overline{w}%
_{m}^{2}-u^{N-2}\overline{w}^{2})dv_{g}+o(1)
\end{equation*}%
and since 
\begin{equation*}
\int_{M}u_{m}^{N-2}\overline{w}_{m}^{2}dv_{g}=1
\end{equation*}%
we obtain%
\begin{equation*}
\left( 1-2^{-\frac{4}{n}}\left( K_{2}^{2}+\epsilon \right) \mu
_{2}(M,g)\right) \int_{M}u^{N-2}\overline{w}^{2}dv_{g}\geq 1-2^{-\frac{4}{n}%
}\left( K_{2}^{2}+\epsilon \right) \mu _{2}(M,g)+o(1)\text{.}
\end{equation*}%
So if 
\begin{equation*}
K_{2}^{2}\mu _{2}(M,g)2^{-\frac{4}{n}}<1
\end{equation*}%
we choose $\epsilon >0$ sufficiently small and get%
\begin{equation*}
\int_{M}u^{N-2}\overline{w}^{2}dv_{g}\geq 1\text{.}
\end{equation*}%
The inequality(\ref{23'}), the Lieb-Brezis lemma (\cite{3}) and the strong
convergence of the sequence $(\overline{w}_{m})_{m}$ to $\overline{w}$ in $%
H_{q}^{2}(M)$ , $q=0,1$, we get 
\begin{equation*}
\int_{M}u_{m}^{N-2}(\overline{w}_{m}-\overline{w})^{2}dv_{g}\leq 2^{-\frac{4%
}{n}}\left( K_{2}^{2}+\epsilon \right) \left( \left\Vert \Delta \overline{w}%
_{m}\right\Vert _{2}^{2}-\left\Vert \Delta \overline{w}\right\Vert
_{2}^{2}\right) +o(1)
\end{equation*}%
\begin{equation*}
\leq 2^{-\frac{4}{n}}\left( K_{2}^{2}+\epsilon \right) \left( \lambda
_{2,m}\int_{M}u_{m}^{N-2}\overline{w}_{m}^{2}dv_{g}-\mu
_{2}(M,g)\int_{M}u^{N-2}\overline{w}^{2}dv_{g}\right) 
\end{equation*}%
\begin{equation*}
+o(1)
\end{equation*}%
\begin{equation*}
\leq 2^{-\frac{4}{n}}\left( K_{2}^{2}+\epsilon \right) \left[ \left( \lambda
_{2,m}-\mu _{2}(M,g)\right) \int_{M}u_{m}^{N-2}\overline{w}%
_{m}^{2}dv_{g}\right. 
\end{equation*}%
\begin{equation*}
+\left. \mu _{2}(M,g)\int_{M}\left( u_{m}^{N-2}\overline{w}_{m}^{2}-u^{N-2}%
\overline{w}^{2}\right) dv_{g}\right] +o(1)\text{.}
\end{equation*}%
Since $\lambda _{2,m}\rightarrow \mu _{2}(M,g)$ as $m\rightarrow +\infty $
and 
\begin{equation}
\int_{M}\left( u_{m}^{N-2}\overline{w}_{m}^{2}-u^{N-2}\overline{w}%
^{2}\right) dv_{g}=\int_{M}u_{m}^{N-2}\left( \overline{w}_{m}-\overline{w}%
\right) ^{2}dv_{g}+o(1)  \label{24}
\end{equation}%
we obtain 
\begin{equation*}
\left( 1-2^{-\frac{4}{n}}\left( K_{2}^{2}+\epsilon \right) \mu
_{2}(M,g)\right) \int_{M}u_{m}^{N-2}\left( \overline{w}_{m}-\overline{w}%
\right) ^{2}dv_{g}\leq o(1)\text{.}
\end{equation*}%
So if 
\begin{equation*}
K_{2}^{2}\mu _{2}(M,g)2^{-\frac{4}{n}}<1
\end{equation*}%
we get%
\begin{equation*}
\lim_{m\rightarrow +\infty }\int_{M}u_{m}^{N-2}\left( \overline{w}_{m}-%
\overline{w}\right) ^{2}dv_{g}=0\text{.}
\end{equation*}%
Hence by the equality(\ref{24}), we get 
\begin{equation*}
\int_{M}u_{m}^{N-2}\overline{w}_{m}^{2}dv_{g}\rightarrow \int_{M}u^{N-2}%
\overline{w}^{2}dv_{g}\text{.}
\end{equation*}%
So 
\begin{equation*}
\int_{M}u^{N-2}\overline{w}^{2}dv_{g}=0
\end{equation*}%
and since 
\begin{equation*}
\int_{M}u^{N-2}v^{2}dv_{g}=\int_{M}u^{N-2}w^{2}dv_{g}=1
\end{equation*}%
and 
\begin{equation*}
\lambda ^{2}+\mu ^{2}=1,\lambda \mu \neq 0
\end{equation*}%
it follows that%
\begin{equation*}
\int_{M}u^{N-2}vwdv_{g}=0\text{.}
\end{equation*}%
Thus the functions $u^{\frac{N-2}{2}}v$, $u^{\frac{N-2}{2}}w$ are linearly
independent.
\end{proof}

\section{An estimation to $\protect\mu _{2}(M,g)$}

Mimicking which is done in \cite{1}, we establish the following lemma.

\begin{lemma}
\label{lem3} If the manifold $(M,g)$ is of dimensional $n\geq 12$, then $\mu
_{2}(M,g)<\left[ \mu _{1}(M,g)^{\frac{n}{4}}+\left( K_{2}^{-2}\right) ^{%
\frac{n}{4}}\right] ^{\frac{4}{n}}$.
\end{lemma}

To prove this lemma, we need the following elementary inequality.

\begin{lemma}
\label{lem4} \cite{1} For any real numbers $x>0$, $y>0$ and $p>2$, there is
a constant $C>0$ such that 
\begin{equation*}
(x+y)^{p}\leq x^{p}+y^{p}+C(x^{p-1}y+xy^{p-1})\text{.}
\end{equation*}
\end{lemma}

\begin{proof}
(of Lemma\cite{lem3}) Let $x_{o}\in M$ , $\delta >0$ sufficiently small and $%
B_{x_{o}}(\delta )$ the ball of center $x_{o}$ and of radius $\delta $ and $%
\eta $ a $C^{\infty }$-function 
\begin{equation*}
\eta (x)=\left\{ 
\begin{array}{c}
1\text{ \ \ if }x\in B_{x_{o}}(\delta ) \\ 
0\text{ \ if }x\notin B_{x_{o}}(2\delta )\text{\ }%
\end{array}%
\right. \text{.}
\end{equation*}

Put 
\begin{equation*}
\varphi _{\epsilon }=\eta (r^{2}+\epsilon ^{2})^{-\frac{n-4}{2}}
\end{equation*}%
where $\eta $ is a bumping function, obviously $\varphi _{\epsilon }\in
H_{2}^{2}(M)$.

For any $n>6$ and $\epsilon $ $\rightarrow 0$, a calculation done in \cite{4}
leads to 
\begin{equation}
Y(\varphi _{\epsilon })\rightarrow K^{-2}-\epsilon ^{2}c_{o}+O(\epsilon ^{3})
\label{26}
\end{equation}%
where%
\begin{equation*}
Y(v)=\frac{\int_{M}vP_{g}(v)dv_{g}}{\left( \int_{M}\left\vert v\right\vert
^{N}dv_{g}\right) ^{\frac{2}{N}}}
\end{equation*}%
$c_{o}>0$

and 
\begin{equation*}
K_{2}^{-2}=\frac{n(n+2)(n-2)(n-4)}{16}\omega _{n-1}^{\frac{4}{n}}\text{.}
\end{equation*}%
$\omega _{n-1}$ denotes the volume of the unit Euclidean sphere.

Consider the function 
\begin{equation}
v_{\epsilon }=c_{\epsilon }\varphi _{\epsilon }  \label{27}
\end{equation}%
with $c_{\epsilon }>0$ is such that $\int_{M}v_{\epsilon }^{N}dv_{g}=1$ .
Standard computations give 
\begin{equation}
c_{\epsilon }=c_{o}\epsilon ^{\frac{n-4}{2}}  \label{28}
\end{equation}%
with $c_{o}>0.$

Denote also by $v$ a smooth positive solution of the equation 
\begin{equation*}
P_{g}(v)=\mu _{1}(M,g)v^{N-1}
\end{equation*}

with $\left\Vert v\right\Vert _{N}=1$. .

Put 
\begin{equation*}
u_{\epsilon }=Y(v_{\epsilon })^{\frac{1}{N-2}}v_{\epsilon }+\mu _{1}(M,g)^{%
\frac{1}{N-2}}v\text{.}
\end{equation*}%
For any $(\lambda ,\mu )\in \mathbb{R}^{2}-\left\{ \left( 0,0\right)
\right\} $, we have 
\begin{equation*}
\int_{M}(\lambda v_{\epsilon }+\mu v)P_{g}(\lambda v_{\epsilon }+\mu
v)dv_{g})^{2}dv_{g}
\end{equation*}%
\begin{equation*}
=\int_{M}\left( \lambda ^{2}v_{\epsilon }P(v_{\epsilon })dv_{g}+\mu
^{2}vP(v)+2\lambda \mu v_{\epsilon }P(v)\right) dv_{g}\text{.}
\end{equation*}%
Since $\int_{M}vP_{g}(v)dv_{g}=\mu _{1}(M,g)$ and $\int_{M}v_{\epsilon
}P_{g}(v_{\epsilon })dv_{g}=Y(v_{\epsilon })$, we get 
\begin{equation*}
\int_{M}(\lambda v_{\epsilon }+\mu v)P_{g}(\lambda v_{\epsilon }+\mu
v)dv_{g}=\lambda ^{2}Y(v_{\epsilon })+\mu ^{2}\mu _{1}(M,g)
\end{equation*}%
\begin{equation*}
+2\lambda \mu \mu _{1}(M,g)\int_{M}v_{\epsilon }v^{N-1}dv_{g}\text{.}
\end{equation*}%
and

\begin{equation*}
\int_{M}u_{\epsilon }^{N-2}(\lambda v_{\epsilon }+\mu v)^{2}dv_{g}=\lambda
^{2}\int_{M}u_{\epsilon }^{N-2}v_{\epsilon }^{2}dv_{g}+\mu
^{2}\int_{M}u_{\epsilon }^{N-2}v^{2}dv_{g}
\end{equation*}%
\begin{equation*}
+2\lambda \mu \int_{M}u_{\epsilon }^{N-2}v_{\epsilon }vdv_{g}
\end{equation*}%
\begin{equation*}
\geq \lambda ^{2}Y(v_{\epsilon })\int_{M}v_{\epsilon }^{N}dv_{g}+\mu ^{2}\mu
_{1}(M,g)\int_{M}v^{N}dv_{g}+2\lambda \mu \int_{M}u_{\epsilon
}^{N-2}v_{\epsilon }vdv_{g}
\end{equation*}%
\begin{equation*}
=\lambda ^{2}Y(v_{\epsilon })+\mu ^{2}\mu _{1}(M,g)+2\lambda \mu
\int_{M}u_{\epsilon }^{N-2}v_{\epsilon }vdv_{g}\text{.}
\end{equation*}%
We have also 
\begin{equation*}
\int_{M}u_{\epsilon }^{N-2}v_{\epsilon }vdv_{g}\geq \mu
_{1}(M,g)\int_{M}v_{\epsilon }v^{N-1}dv_{g}
\end{equation*}%
so, if $\lambda \mu \geq 0$ 
\begin{equation*}
\frac{\int_{M}(\lambda v_{\epsilon }+\mu v)P_{g}(\lambda v_{\epsilon }+\mu
v)dv_{g}}{\int_{M}u_{\epsilon }^{N-2}(\lambda v+\mu v_{\epsilon })^{2}dv_{g}}%
\leq 1\text{.}
\end{equation*}

In the case $\lambda \mu <0$ and $N-2\in (0,1]$ i.e. $n\geq 12$, we have%
\begin{equation*}
u_{\epsilon }^{N-2}=\left( Y(v_{\epsilon })^{\frac{1}{N-2}}v_{\epsilon }+\mu
_{1}(M,g)^{\frac{1}{N-2}}v\right) ^{N-2}
\end{equation*}%
\begin{equation*}
\leq Y(v_{\epsilon })v_{\epsilon }^{N-2}+\mu _{1}(M,g)v^{N-2}\text{.}
\end{equation*}%
Consequently%
\begin{equation*}
\int_{M}u_{\epsilon }^{N-2}(\lambda v_{\epsilon }+\mu v)^{2}dv_{g}=\lambda
^{2}\int_{M}u_{\epsilon }^{N-2}v_{\epsilon }^{2}dv_{g}+\mu
^{2}\int_{M}u_{\epsilon }^{N-2}v^{2}dv_{g}
\end{equation*}%
\begin{equation*}
+2\lambda \mu \int_{M}u_{\epsilon }^{N-2}v_{\epsilon }vdv_{g}\geq \lambda
^{2}Y(v_{\epsilon })+\mu ^{2}\mu _{1}(M,g)
\end{equation*}%
\begin{equation*}
+2\lambda \mu Y(v_{\epsilon })\int_{M}v_{\epsilon }^{N-1}vdv_{g}+2\lambda
\mu \mu _{1}(M,g)\int_{M}v_{\epsilon }v^{N-1}dv_{g}
\end{equation*}%
\begin{equation*}
\geq \lambda ^{2}Y(v_{\epsilon })+\mu ^{2}\mu _{1}(M,g)-C\left(
\int_{M}v_{\epsilon }^{N-1}vdv_{g}+\int_{M}v_{\epsilon }v^{N-1}dv_{g}\right)
\end{equation*}%
where $C>0$ is a constant independent of $\epsilon $.

Now taking account of (\ref{27}) and (\ref{28}) we get, for any $(\lambda
,\mu )\in \mathbb{R}^{2}-\left\{ \left( 0,0\right) \right\} $,%
\begin{equation*}
\frac{\int_{M}(\lambda v_{\epsilon }+\mu v)P_{g}(\lambda v_{\epsilon }+\mu
v)dv_{g}}{\int_{M}u_{\epsilon }^{N-2}(\lambda v+\mu v_{\epsilon })^{2}dv_{g}}%
\leq 1+O(\epsilon ^{\frac{n-4}{2}})\text{.}
\end{equation*}

By Lemma\ref{lem4}, we obtain%
\begin{equation*}
\int_{M}u_{\epsilon }^{N}dv_{g}\leq Y(v_{\epsilon })^{\frac{n}{4}%
}\int_{M}v_{\epsilon }^{N}dv_{g}+\mu _{1}(M,g)^{\frac{n}{4}%
}\int_{M}v^{N}dv_{g}\geq 
\end{equation*}%
\begin{equation*}
+C\left( \int_{M}v_{\epsilon }^{N-1}vdv_{g}+\int_{M}vv_{\epsilon
}^{N-1}dv_{g}\right) 
\end{equation*}%
\begin{equation*}
=Y(v_{\epsilon })^{\frac{n}{4}}+\mu _{1}(M,g)^{\frac{n}{4}}+C\left(
\int_{M}v_{\epsilon }^{N-1}vdv_{g}+\int_{M}v_{\epsilon }v^{N-1}dv_{g}\right) 
\text{.}
\end{equation*}%
And by the relation (\ref{26}), we deduce that%
\begin{equation*}
\left( \int_{M}u_{\epsilon }^{N}dv_{g}\right) ^{\frac{4}{n}}\leq \left[ \mu
_{1}(M,g)^{\frac{n}{4}}+(K_{2}^{-2})^{\frac{n}{4}}-c.\epsilon
^{2}+o(\epsilon ^{3})+o(\epsilon ^{\frac{n-4}{2}})\right] ^{\frac{4}{n}}
\end{equation*}%
\begin{equation*}
=\left[ \mu _{1}(M,g)^{\frac{n}{4}}+(K_{2}^{-2})^{\frac{n}{4}}\right] ^{%
\frac{4}{n}}-c.\epsilon ^{2}+o(\epsilon ^{3})\text{ \ \ \ \ \ \ (}\frac{n-4}{%
2}\geq 4\text{).}
\end{equation*}%
where $c>0$ is a constant.

Hence, for any $(\lambda ,\mu )\in \mathbb{R}^{2}-\left\{ \left( 0,0\right)
\right\} $,%
\begin{equation*}
\frac{\int_{M}(\lambda v_{\epsilon }+\mu v)P_{g}(\lambda v_{\epsilon }+\mu
v)dv_{g}}{\int_{M}u_{\epsilon }^{N-2}(\lambda v+\mu v_{\epsilon })^{2}dv_{g}}%
\left( \int_{M}u_{\epsilon }^{N}dv_{g}\right) ^{\frac{4}{n}}\leq
\end{equation*}%
\begin{equation*}
\left[ \mu _{1}(M,g)^{\frac{n}{4}}+(K_{2}^{-2})^{\frac{n}{4}}\right] ^{\frac{%
4}{n}}-c.\epsilon ^{2}+O(\epsilon ^{3})
\end{equation*}%
so

\begin{equation*}
\mu _{2}(M,g)<\left[ \mu _{1}(M,g)^{\frac{n}{4}}+\left( K_{2}^{-2}\right) ^{%
\frac{n}{4}}\right] ^{\frac{4}{n}}\text{.}
\end{equation*}
\end{proof}

\section{Strong convergence}

\begin{lemma}
\label{lem5} Suppose that $\mu _{2}K_{2}^{2}2^{-\frac{4}{n}}<1.$ Then the
sequence $\left( v_{m}\right) _{m}$ (resp. $\left( w_{m}\right) _{m}$ ) of
solutions of the equations(\ref{17}) (resp. of solutions of the equations (%
\ref{18})) has a bounded subsequence on $M$.
\end{lemma}

\begin{proof}
Let as in the section2 $u_{m}\in L_{+}^{N}(M)$ and $v_{m},w_{m}\in
H_{2}^{2}(M)$ solutions respectively of the equations 
\begin{equation}
P(v_{m})=\lambda _{1,m}u_{m}^{N-2}v_{m}  \label{29}
\end{equation}%
and 
\begin{equation}
P(w_{m})=\lambda _{2,m}u_{m}^{N-2}w_{m}\text{.}  \label{29'}
\end{equation}%
Set%
\begin{equation*}
S_{m}=\left\{ \lambda _{m}v_{m}+\mu _{m}w_{m}\text{:\ }\left( \lambda
_{m},\mu _{m}\right) \in R^{2}\text{, }\lambda _{m}^{2}+\mu _{m}^{2}=1\text{%
, }\lambda _{m}\mu _{m}>\alpha >0\text{\ }\right\}
\end{equation*}%
\begin{equation*}
S=\left\{ \lambda v+\mu w\text{:\ }\left( \lambda ,\mu \right) \in R^{2}%
\text{, }\lambda ^{2}+\mu ^{2}=1\right\}
\end{equation*}%
and let $\overline{w}_{m}=\lambda _{m}v_{m}+\mu _{m}w_{m}$ , $\overline{w}%
=\lambda v+\mu w$ where up to a subsequence $\left( \lambda _{m},\mu
_{m}\right) \rightarrow \left( \lambda ,\mu \right) $.

First we are going to show that the sequence $\left( \overline{w}_{m}\right)
_{m}$ is uniformally bounded on the manifold $M.$ Suppose by contradiction
that $(\overline{W}_{m})_{m}$ is unbounded. Then, for every $m$ there exists
a point $x_{m}\in M$ such that 
\begin{equation*}
\overline{w}_{m}(x_{m})=\max_{x\in M}\overline{w}_{m}=\xi _{m}\rightarrow
+\infty \text{ as }m\rightarrow +\infty \text{.}
\end{equation*}%
Given $\delta >0$ less than the injectivity radius of $(M,g)$, we let $%
\widetilde{w}_{m}$ and $\widetilde{u}_{m}$ be the functions defined on the
Euclidean ball of center $0$ and radius $\delta \xi _{m}$, $B_{o}(\delta \xi
_{m})$, \ by%
\begin{equation*}
\widetilde{w}_{m}(x)=\frac{1}{\xi _{m}}\overline{w}_{m}(\exp _{x_{m}}(\frac{x%
}{\xi _{m}}))
\end{equation*}%
where $\exp _{m}$ is the exponential map at $x_{m}$. Denote by 
\begin{equation*}
g_{m}(x)=\left( \exp _{x_{m}}\right) ^{\ast }g(\frac{x}{\xi _{m}})
\end{equation*}%
the Riemannian metric on the ball $B_{o}(\delta \xi _{m})$. Clearly, if $E$
is the Euclidean metric, \ $g_{m}\rightarrow E$ in $C^{2}$ on any compact
set.

Now, since the functions $v_{m}$ and $w_{m}$ are solutions respectively of
the equations (\ref{29}) and (\ref{29'}) then multiplying by $\overline{w}%
_{m}$ $\in H_{2}^{2}(M)$ and integrating over the geodesic ball $\exp
_{x_{m}}(B_{o}(\delta \xi _{m}))$, we obtain%
\begin{equation*}
\int_{B_{o}(\delta \xi _{m})}\widetilde{w}_{m}\Delta _{g_{m}}^{2}\widetilde{w%
}_{m}dv_{g_{m}}+\frac{\alpha }{\xi _{m}^{2}}\int_{B_{o}(\delta \xi _{m})}%
\widetilde{w}_{m}\Delta \widetilde{w}_{m}dv_{g_{m}}
\end{equation*}%
\begin{equation*}
+\frac{a}{\xi _{m}^{4}}\int_{B_{o}(\delta \xi _{m})}\widetilde{w}%
_{m}^{2}dv_{g_{m}}\leq \lambda _{2,m}\int_{B_{o}(\delta \xi _{m})}\widetilde{%
u}_{m}^{N-2}\widetilde{w}_{m}^{2}dv_{g_{m}}\text{.}
\end{equation*}%
Thus the fact that $\left\vert \widetilde{w}_{m}\right\vert \leq 1$ on $%
B_{o}(\delta \xi _{m})$ and standard elliptic theory lead after passing to a
subsequence to 
\begin{equation*}
\widetilde{w}_{m}\rightarrow \widetilde{w}\text{ \ \ \ in }C_{loc}^{4}(R^{n})%
\text{.}
\end{equation*}%
Independently, we have, for any $R>0$%
\begin{equation*}
\int_{B_{o}(R\delta )}\widetilde{u}_{m}^{N-2}dx=\int_{B_{o}(R\delta
)}u_{m}^{N-2}(\exp _{x_{m}}(\frac{x}{R}))\widetilde{w}_{m}^{2}dv_{g_{m}}+o(1)
\end{equation*}%
\begin{equation*}
=\int_{B(x_{m},R\delta )}u_{m}^{N-2}(\exp _{x_{m}}(\frac{x}{R}))dv_{g}+o(1)
\end{equation*}%
\begin{equation*}
\leq \int_{M}u_{m}^{N-2}(x)dv_{g}+o(1)\text{.}
\end{equation*}%
So $\widetilde{u}_{m}\rightarrow \widetilde{u}$ weakly in $L_{loc}^{N}(R^{n})
$ 
\begin{equation*}
\int_{B_{o}(R\delta )}u_{m}^{N}(\exp _{x_{m}}(\frac{x}{R}))dv_{g_{m}}=%
\int_{B(x_{m},R\delta )}u_{m}^{N}(x)dv_{g}\leq \int_{M}u_{m}^{N}(x)dv_{g}
\end{equation*}%
Now letting $m\rightarrow \infty $, we get%
\begin{equation*}
\int_{R^{n}}\left( \Delta _{E}\widetilde{w}\right) ^{2}dx\leq \mu
_{2}\int_{R^{n}}u^{N-2}\widetilde{w}^{2}dx
\end{equation*}%
and by the Sobolev inequality given by Corollary\ref{cor1}, we obtain that 
\begin{equation*}
\int_{R^{n}}\left( \Delta _{E}\widetilde{w}\right) ^{2}dx\leq \mu _{2}2^{-%
\frac{4}{n}}K_{2}^{2}\int_{R^{n}}\left( \Delta _{E}\widetilde{w}\right)
^{2}dx\left( \int_{R^{n}}u^{N}dx\right) ^{\frac{2}{N}}
\end{equation*}%
\begin{equation*}
\leq \mu _{2}2^{-\frac{4}{n}}K_{2}^{2}\int_{R^{n}}\left( \Delta _{E}%
\widetilde{w}\right) ^{2}dx\text{.}
\end{equation*}%
Consequently 
\begin{equation*}
\mu _{2}K_{2}^{2}2^{-\frac{4}{n}}\geq 1
\end{equation*}%
which contradicts the inequalityof the hypothesis, so the sequence $\left( 
\overline{w}_{m}\right) $ is bounded on $M$.
\end{proof}

\begin{corollary}
\label{co1}The sequence $\left( v_{m}\right) _{m}$ (resp. $(w_{m})_{m}$
)given in (\ref{17}) (resp. in (\ref{18})) converges strongly in $L^{N}(M)$.
\end{corollary}

\begin{proof}
Let $\epsilon >0$, the H\"{o}lder inequality leads to%
\begin{equation*}
\int_{M}\left\vert v_{n}-v\right\vert ^{N}dv_{g}\leq \left(
\int_{M}\left\vert v_{n}-v\right\vert ^{N-\epsilon }dv_{g}\right) ^{\frac{1}{%
2}}\left( \int_{M}\left\vert v_{n}+v\right\vert ^{N+\epsilon }dv_{g}\right)
^{\frac{1}{2}}\text{.}
\end{equation*}%
By the boundedness of the sequence $(v_{n})_{n}$ in $M$ and the strong
convergence of the latter to $v$ in $L^{N-\epsilon }(M)$, we get that $%
v_{n}\rightarrow v$ in $L^{N}(M)$. The same is also true for the sequence $%
(w_{n})_{n}$.
\end{proof}

\begin{corollary}
The functions $u^{\frac{N-2}{2}}v$ and $u^{\frac{N-2}{2}}w$ are linearly
independent.
\end{corollary}

Indeed, since the sequence $(v_{m})_{m}$ ( resp.$(v_{m})_{m}$ ) converges
strongly to $v$ (resp.to $w$) in $L^{N}(M)$, we pass to the limit in the
last equality in (\ref{19}) and get $\int_{M}u^{N-2}vwdv_{g}=0$

As a corollary of Lemma \ref{lem5} and Corollary \ref{co1}, we obtain our
main result

\begin{theorem}
If the Einsteinian manifold $(M,g)$ is of dimension $n\geq 12$, then $\mu
_{2}(M,g)$ is attained by a generalized metric.
\end{theorem}

\section{Nodals solutions}

The same arguments as in the proof of Lemma3.3 \cite{1} \ allow us to state.

\begin{lemma}
\label{lem6} Let $u\in L_{+}^{N}(M)$ with $\left\Vert u\right\Vert _{N}=1.$%
Suppose that $w_{1}$, $w_{2}\in H_{2}^{2}(M)-\left\{ 0\right\} $, such that $%
w_{1}\geq 0$, $w_{2}\geq 0$ satisfy 
\begin{equation}
\int_{M}w_{1}P(w_{1})dv_{g}\leq \mu _{2}(M,g)\int_{M}u^{N-2}w_{1}^{2}dv_{g}
\label{30}
\end{equation}%
\begin{equation}
\int_{M}w_{2}P(w_{2})dv_{g}\leq \mu _{2}(M,g)\int_{M}u^{N-2}w_{2}^{2}dv_{g}%
\text{.}  \label{31}
\end{equation}%
If $(M-w_{1}^{-1}(0))\cap (M-w_{2}^{-1}(0))$ has measure $0$, then there
exist constants $a>0$ and $b>0$ such that $u=aw_{1}+bw_{2}$ and the
equalities in (\ref{30}) and (\ref{31}) hold.
\end{lemma}

Now we establish the existence of a nodal solution.

\begin{theorem}
Let $v$ and $w$ as in the Proposition(\ref{p3}) and suppose that the scalar
curvature of $(M;g)$ is positive and $\mu _{2}(M,g)\neq 0$ and attained by a
general metric $\widetilde{g}=u^{N-2}g$ with $u\in L_{+}^{N}(M)$. Then $%
u=\left\vert w\right\vert $ \ and in particular the equation 
\begin{equation}
P_{g}(w)=\mu _{2}(M,g)\left\vert w\right\vert ^{N-2}w  \label{32}
\end{equation}%
has a nodal solution.
\end{theorem}

\begin{proof}
Let $v$ and $w$ as in the Proposition(\ref{p3}). Without lost of generality,
we choose $u\in L_{+}^{N}(M)$ with $\int_{M}u^{N}dv_{g}=1$, hence $\lambda
_{2}(\widetilde{g})=\mu _{2}(M,g).$ As in the proof of Proposition\ref{pro1}
we have $\lambda _{1}(\widetilde{g})<\lambda _{2}(\widetilde{g})$. Suppose
that the solution of the equation(\ref{32}) is not nodal, by taking $-$ $w$
if $w$ is non positive, we assume that $w\geq 0$. On the other hand since
the scalar curvature of $(M,g)$ is positive, by Proposition(\ref{p2}) the
equation 
\begin{equation*}
P_{g}(v)=\lambda _{1}(\widetilde{g})u^{N-2}v
\end{equation*}%
has a positive solution and by Proposition\ref{p3} \ with the constraints 
\begin{equation*}
\int_{M}u^{N-2}v^{2}dv_{g}=\int_{M}u^{N-2}w^{2}dv_{g}=1
\end{equation*}%
and 
\begin{equation*}
\int_{M}u^{N-2}vwdv_{g}=0\text{.}
\end{equation*}%
This latter equality implies that the set $(M-v^{-1}(0))\cap (M-w^{-1}(0))$
has measure $0$. So by Lemma\ref{lem6}, we get equalities in (\thinspace \ref%
{32}), a contradiction with the fact that $\lambda _{1}(\widetilde{g}%
)<\lambda _{2}(\widetilde{g})$. Consequently $w$ is a nodal function.

Suppose that the compact manifold $M$ splits into two non empty disjoint
domains $\Omega _{1}$and $\Omega _{2}$ such that $M=\Omega _{1}\cup \Omega
_{2}\cup \digamma $ \ with measure$(\digamma )=0$. Let $v_{1}$ and $v_{2}$
be positive solutions to the equation $P_{g}(v_{i})=\lambda _{2}u^{N-2}v_{i}$%
, such that $v_{i}=0$ and $\Delta v_{i}=0$ on $\partial \Omega _{i}$, where $%
\lambda _{2}$ is the second eigenvalue of the Paneitz-Branson operator $P_{g}
$. By Lemma\ref{lem6}, there exist constants $a>0$ and $b>0$ such that $%
u=av_{1}+bv_{2}$. It follows that $u$ is of class $C^{o,\alpha }(M)$ with $%
\alpha \in \left( 0,N-2\right) $. Observe that the nodal set $%
u^{-1}(0)\subset v_{1}^{-1}(0)\cap v_{2}^{-1}(0)\subset \digamma .$

Now we follows the proof in \cite{1}. Let $h\in C^{\infty }(M)$ with support
in $M-u^{-1}(0)$ and put $u_{t}=u+th$. Since $u$ is continuous and positive
on the support of $h$, then $u_{t}>0$ for $t$ close to $0$. The same
arguments as the proof in the Proposition3.3 in \cite{1} we obtain that $%
\left\vert w\right\vert =u$ on $M-u^{-1}(0)$. Independently since the nodal
set $u^{-1}(0)$ is negligible and $u$, $\left\vert w\right\vert $ are
continuous, then $\left\vert w\right\vert =u$ on $M$.
\end{proof}

\begin{corollary}
$\mu _{2}(M,g)$ is not attained by a classical conformal metric.
\end{corollary}

Since if it is not the case, $u>0$ and $w$ such that $\left\vert
w\right\vert =u$ is not nodal.

\end{document}